\documentclass{article}
\usepackage{amssymb}
\usepackage{amsmath}

\newtheorem{theorem}{Theorem}

\newtheorem{definition}{Definition}

\newtheorem{lemma}{Lemma}

\newtheorem{proposition}{Proposition}
\newtheorem{remark}{Remark}

\begin{document}

\title{\textbf{Strong approximation of almost periodic functions }}
\author{\textbf{W\l odzimierz \L enski \ and Bogdan Szal} \\
University of Zielona G\'{o}ra\\
Faculty of Mathematics, Computer Science and Econometrics\\
65-516 Zielona G\'{o}ra, ul. Szafrana 4a\\
P O L A N D\\
W.Lenski@wmie.uz.zgora.pl ,B.Szal @wmie.uz.zgora.pl}
\date{}
\maketitle

\begin{abstract}
We consider summability methods generated by the class $GM\left( _{2}\beta
\right) .$ We generalize some related results of \ P. Pych-Taberska [Studia
Math. XCVI (1990),91-103] on strong approximation of almost periodic
functions\ by their Fourier series and S. M. Mazhar and V.Totik [J. Approx.
Theory, 60(1990), 174--182] on approximation of periodic functions by matrix
means of their Fourier series .

\ \ \ \ \ \ \ \ \ \ \ \ \ \ \ \ \ \ \ \ 


\ \ \ \ \ \ \ \ \ \ \ \ \ \ \ \ \ \ \ 

\end{abstract}

\section{Introduction}

Let $S^{p}\;\left( 1<p\leq \infty \right) $ be the class of all almost
periodic functions in the sense of Stepanov with the norm%
\begin{equation*}
\Vert f\Vert _{S^{p}}:=\left\{ 
\begin{array}{c}
\sup\limits_{u}\left\{ \frac{1}{\pi }\int_{u}^{u+\pi }\mid f(t)\mid
^{p}dt\right\} ^{1/p}\text{ \ \ when \ \ }1<p<\infty \\ 
\sup\limits_{u}\mid f(u)\mid \text{ \ \ when \ \ }p=\infty .%
\end{array}%
\right.
\end{equation*}%
Denote yet by $C_{2\pi }$ the class of all $2\pi $--periodic functions
continuous over $Q=$ $[-\pi ,\pi ]$ with the norm%
\begin{equation*}
\Vert f\Vert _{C_{2\pi }}:=\sup\limits_{t\in Q}\mid f(t)\mid .
\end{equation*}%
Suppose that the Fourier series of $f\in S^{p}$ has the form 
\begin{equation*}
Sf\left( x\right) =\sum_{\nu =-\infty }^{\infty }A_{\nu }\left( f\right)
e^{i\lambda _{\nu }x},\text{ \ \ where \ }A_{\nu }\left( f\right)
=\lim_{L\rightarrow \infty }\frac{1}{L}\int_{0}^{L}f(t)e^{-i\lambda _{\nu
}t}dt,
\end{equation*}%
with the partial sums\ 
\begin{equation*}
S_{\gamma _{k}}f\left( x\right) =\sum_{\left\vert \lambda _{\nu }\right\vert
\leq \gamma _{k}}A_{\nu }\left( f\right) e^{i\lambda _{\nu }x}
\end{equation*}%
and that $0=\lambda _{0}<\lambda _{\nu }<\lambda _{\nu +1}$ if $\nu \in 
\mathbb{N}
=\left\{ 1,2,3...\right\} ,$ $\lim \lambda _{\nu }=\infty ,$ $\lambda _{-\nu
}=-\lambda _{\nu }$\ $\left\vert A_{\nu }\right\vert +\left\vert A_{-\nu
}\right\vert >0.$ Let\ $\Omega _{\alpha ,p}$ , with some fixed positive $%
\alpha $ , be the set of functions of class $S^{p}$ bounded on $%
\mathbb{R}
=\left( -\infty ,\infty \right) $ whose Fourier exponents satisfy the
condition%
\begin{equation*}
\lambda _{\nu +1}-\lambda _{\nu }\geq \alpha \text{ \ \ }\left( \nu \in 
\mathbb{N}
\right) .
\end{equation*}%
In case $f\in \Omega _{\alpha ,p}$\ 
\begin{equation*}
S_{\lambda _{k}}f\left( x\right) =\int_{0}^{\infty }\left\{ f\left(
x+t\right) +f\left( x-t\right) \right\} \Psi _{\lambda _{k},\lambda
_{k}+\alpha }\left( t\right) dt,
\end{equation*}%
where%
\begin{equation*}
\Psi _{\lambda ,\eta }\left( t\right) =\frac{2\sin \frac{\left( \eta
-\lambda \right) t}{2}\sin \frac{\left( \eta +\lambda \right) t}{2}}{\pi
\left( \eta -\lambda \right) t^{2}}\text{ \ \ }\left( 0<\lambda <\eta ,\text{
\ }\left\vert t\right\vert >0\right) .
\end{equation*}

Let $A:=\left( a_{n,k}\right) $ be an infinite matrix of real nonnegative
numbers such that%
\begin{equation}
\sum_{k=0}^{\infty }a_{n,k}=1\text{, where }n=0,1,2,...\text{\ ,}  \label{T1}
\end{equation}%
and let the\ $A-$transformation$\ $of \ $\left( S_{\gamma _{k}}f\right) $ be
given by%
\begin{equation*}
T_{n,A,\gamma }^{\text{ }}f\left( x\right) :=\sum_{k=0}^{\infty
}a_{n,k}S_{\gamma _{k}}f\left( x\right) \text{ \ \ \ }\left( \text{ \ }%
n=0,1,2,...\right) .
\end{equation*}

Let us consider the strong mean 
\begin{equation}
T_{n,A,\gamma }^{q}f\left( x\right) =\left\{ \sum_{k=0}^{\infty
}a_{n,k}\left\vert S_{\gamma _{k}}f\left( x\right) -f\left( x\right)
\right\vert ^{q}\right\} ^{1/q}\text{ \ \ \ }\left( q>0\right) \text{.}
\label{S1}
\end{equation}%
If $f\in C_{2\pi },$\ then as usually 
\begin{equation*}
Sf\left( x\right) =\frac{a_{o}(f)}{2}+\sum_{\nu =1}^{\infty }(a_{\nu
}(f)\cos kx+b_{\nu }(f)\sin kx)
\end{equation*}%
and instead of \ $S_{\gamma _{k}}f$\ we will consider the partial sums\ 
\begin{equation*}
S_{k}f\left( x\right) =\frac{a_{o}(f)}{2}+\sum_{\nu =1}^{k}(a_{\nu }(f)\cos
kx+b_{\nu }(f)\sin kx).
\end{equation*}%
Thus, instead of $T_{n,A,\gamma }^{\text{ }}f$\ and $T_{n,A,\gamma }^{q}f$\
we will consider the quantities $T_{n,A}^{\text{ }}f$\ and $T_{n,A}^{q}f$
defined by the formulas\ 
\begin{equation}
T_{n,A}^{\text{ }}f\left( x\right) :=\sum_{k=0}^{\infty }a_{n,k}S_{k}f\left(
x\right) \text{ \ \ \ }\left( \text{ \ }n=0,1,2,...\right)  \label{S2}
\end{equation}%
and%
\begin{equation}
T_{n,A}^{q}f\left( x\right) =\left\{ \sum_{k=0}^{\infty }a_{n,k}\left\vert
S_{k}f\left( x\right) -f\left( x\right) \right\vert ^{q}\right\} ^{1/q}\text{
\ \ \ }\left( q>0\right) \text{,}  \label{S3}
\end{equation}%
respectively. As measures of approximation by the quantities (\ref{S1}), (%
\ref{S2}) and (\ref{S3}) we use the best approximation of $f$ by
trigonometric polynomials $t_{k}$ of order at most $k$ or by entire
functions $g_{\sigma }$ of exponential type $\sigma $\ bounded on the real
axis, shortly $g_{\sigma }\in B_{\sigma }$ and the modulus of continuity of$%
\ f,$ defined by the formulas%
\begin{equation*}
E_{k}(f)_{C_{2\pi }}=\inf_{t_{k}}\left\Vert f-t_{k}\right\Vert _{C_{2\pi }}
\end{equation*}%
or%
\begin{equation*}
E_{\sigma }(f)_{S^{p}}=\inf_{g_{\sigma }}\left\Vert f-g_{\sigma }\right\Vert
_{S^{p}}
\end{equation*}%
and%
\begin{equation*}
\omega f\left( \delta \right) _{X}=\sup_{\left\vert t\right\vert \leq \delta
}\left\Vert f\left( \cdot +t\right) -f\left( \cdot \right) \right\Vert _{X%
\text{ }},\text{ \ \ }X=C_{2\pi }\text{ or }X=S^{p},
\end{equation*}%
respectively.

In \cite{MT} S. M. Mazhar and V. Totik proved the following theorem:

\begin{theorem}
Let $f\in C_{2\pi }.$ Suppose $A:=\left( a_{n,k}\right) $ satisfies (\ref{T1}%
), $lim_{n\rightarrow \infty }a_{n,0}=0$ and%
\begin{equation*}
a_{n,k}\geq a_{n,k+1}\text{ \ \ }k=0,1,2,...\text{ \ \ }n=0,1,2,...,
\end{equation*}%
then%
\begin{equation*}
\left\Vert T_{n,A}^{\text{ }}f\left( x\right) -f\right\Vert _{C_{2\pi }}\leq
K\sum_{k=0}^{\infty }a_{n,k}\omega f\left( \frac{1}{k+1}\right) _{C_{2\pi }}.
\end{equation*}
\end{theorem}

Recently, L. Leindler \cite{3} defined a new class of sequences named as
sequences of rest bounded variation, briefly denoted by $RBVS$, i.e.,%
\begin{equation}
RBVS=\left\{ a:=\left( a_{n}\right) \in 
\mathbb{C}
:\sum\limits_{k=m}^{\infty }\left\vert a_{k}-a_{k+1}\right\vert \leq
K\left( a\right) \left\vert a_{m}\right\vert \text{ for all }m\in 
\mathbb{N}
\right\} ,  \label{1}
\end{equation}%
where here and throughout the paper $K\left( a\right) $ always indicates a
constant only depending on $a$.

Denote by $MS$ the class of nonincreasing sequences. Then it is obvious that%
\begin{equation*}
MS\subset RBVS.
\end{equation*}%
In \cite{8} L. Leindler considered the class of mean rest bounded variation
sequences $MRBVS$, where%
\begin{equation*}
MRBVS=\left\{ a:=\left( a_{n}\right) \in 
\mathbb{C}
:\right.
\end{equation*}%
\begin{equation}
\left. \sum\limits_{k=m}^{\infty }\left\vert a_{k}-a_{k+1}\right\vert \leq
K\left( a\right) \frac{1}{m}\sum\limits_{k\geq m/2}^{m}\left\vert
a_{k}\right\vert \text{ for all }m\in 
\mathbb{N}
\right\} .  \label{2}
\end{equation}%
It is clear that%
\begin{equation*}
RBVS\subseteq MRBVS.
\end{equation*}%
In \cite{S} the second author proved that $RBVS\neq MRBVS$. Moreover, the
above theorem was generalized for the class $MRBVS$ in \cite{Sz} .

Further, the class of general monotone coefficients, $GM$, is defined as
follows ( see \cite{10}):%
\begin{equation}
GM=\left\{ a:=\left( a_{n}\right) \in 
\mathbb{C}
:\sum\limits_{k=m}^{2m-1}\left\vert a_{k}-a_{k+1}\right\vert \leq K\left(
a\right) \left\vert a_{m}\right\vert \text{ for all }m\in 
\mathbb{N}
\right\} .  \label{3}
\end{equation}%
It is clear%
\begin{equation*}
RBVS\subset GM\text{.}
\end{equation*}%
In \cite{6, 10, 11, 12} was defined the class of $\beta -$general monotone
sequences as follows:

\begin{definition}
Let $\beta :=\left( \beta _{n}\right) $ be a nonnegative sequence. The
sequence of complex numbers $a:=\left( a_{n}\right) $ is said to be $\beta -$%
general monotone, or $a\in GM\left( \beta \right) $, if the relation%
\begin{equation}
\sum\limits_{k=m}^{2m-1}\left\vert a_{k}-a_{k+1}\right\vert \leq K\left(
a\right) \beta _{m}  \label{4}
\end{equation}%
holds for all $m$.
\end{definition}

In the paper \cite{12} Tikhonov considered, among others, the following
examples of the sequences $\beta _{n}:$

(1) $_{1}\beta _{n}=\left\vert a_{n}\right\vert ,$

(2) $_{2}\beta _{n}=\sum\limits_{k=\left[ n/c\right] }^{\left[ cn\right] }%
\frac{\left\vert a_{k}\right\vert }{k}$ for some $c>1$.

It is clear that $GM\left( _{1}\beta \right) =GM.$ Moreover (see \cite[%
Remark 2.1]{12})%
\begin{equation*}
GM\left( _{1}\beta +_{2}\beta \right) \equiv GM\left( _{2}\beta \right) .
\end{equation*}

Consequently, we assume that the sequence $\left( K\left( \alpha _{n}\right)
\right) _{n=0}^{\infty }$ is bounded, that is, that there exists a constant $%
K$ such that%
\begin{equation*}
0\leq K\left( \alpha _{n}\right) \leq K
\end{equation*}%
holds for all $n$, where $K\left( \alpha _{n}\right) $ denote the sequence
of constants appearing in the inequalities (\ref{1})-(\ref{4}) for the
sequences $\alpha _{n}:=\left( a_{nk}\right) _{k=0}^{\infty }$.

Now we can give the conditions to be used later on. We assume that for all $%
n $%
\begin{equation}
\sum\limits_{k=m}^{2m-1}\left\vert a_{n,k}-a_{n,k+1}\right\vert \leq
K\sum\limits_{k=[m/c]}^{[cm]}\frac{a_{n,k}}{k}  \label{5}
\end{equation}%
holds if $\alpha _{n}=\left( a_{n,k}\right) _{k=0}^{\infty }$ belongs to $%
GM\left( _{2}\beta \right) $, for $n=1,2,...$

We have shown in \cite{WLBS} that the following lemma:

\begin{lemma}
(see \cite[Theorem 1.11]{LL} ). Suppose that $n=O\left( r_{n}\right) $ and $%
f\in C_{2\pi }.$ Then, for any number $q>0,$we have%
\begin{equation*}
\left\Vert \left\{ \frac{1}{r_{n}}\sum_{k=n-r_{n}}^{n-1}\left\vert
S_{k}f-f\right\vert ^{q}\right\} ^{1/q}\right\Vert _{C_{2\pi }}\ll
E_{n-r_{n}}(f)_{C_{2\pi }}.
\end{equation*}
\end{lemma}

implies

\begin{theorem}
If $f\in C_{2\pi }$, $q>0$, $\left( a_{n,k}\right) _{k=0}^{\infty }\in
GM\left( _{2}\beta \right) $ for all $n$, (\ref{T1}) and $\lim_{n\rightarrow
\infty }a_{n,0}=0$ hold, then%
\begin{equation}
\left\Vert T_{n,A}^{q}f\right\Vert _{C_{2\pi }}\ll \left\{
\sum_{k=0}^{\infty }a_{n,k}E_{\left[ \frac{k}{2^{\left[ c\right] }}\right]
}^{q}(f)_{C_{2\pi }}\right\} ^{1/q}  \label{6}
\end{equation}%
for some $c>1,$ and%
\begin{equation*}
\left\Vert T_{n,A}^{q}f\right\Vert _{C_{2\pi }}\ll \left\{
\sum_{k=0}^{\infty }a_{n,k}\omega ^{q}f\left( \frac{\pi }{k+1}\right)
_{C_{2\pi }}\right\} ^{1/q},
\end{equation*}%
for $n=0,1,2,...$If we additionally suppose that $\left( a_{n,k}\right)
_{k=0}^{\infty }\in MS,$ then from (\ref{6}) we deduce%
\begin{equation*}
\left\Vert T_{n,A}^{q}f\right\Vert _{C_{2\pi }}\ll \left\{
\sum_{k=0}^{\infty }a_{n,k}E_{k}^{q}(f)_{C_{2\pi }}\right\} ^{1/q},
\end{equation*}%
for $n=0,1,2,...$
\end{theorem}

In this paper we consider the class $GM\left( _{2}\beta \right) $ in
analogical estimate of the quantity $\left\Vert T_{n,A,\gamma
}^{q}f\right\Vert _{S^{p}}.$ Thus we extend the result of S. M. Mazhar and
V. Totik (see \cite[Theorem 1]{MT}) and generalize the following result of
P. Pych-Taberska\ (see \cite[Theorem 5]{PT}):

\begin{theorem}
If $f\in \Omega _{\alpha ,\infty }$ and $q\geq 2,$ then%
\begin{equation*}
\left\Vert T_{n,A,\gamma }^{q}f\right\Vert _{S^{\infty }}\ll \left\{ \frac{1%
}{n+1}\sum_{k=0}^{n}\left[ \omega f\left( \frac{\pi }{k+1}\right)
_{S^{\infty }}\right] ^{q}\right\} ^{1/q}+\frac{\left\Vert f\right\Vert
_{S^{\infty }}}{\left( n+1\right) ^{1/q}},
\end{equation*}%
for $n=0,1,2,...$, where $\gamma =\left( \gamma _{k}\right) $ is a sequence
with $\gamma _{k}=\frac{\alpha k}{2},$ $a_{n,k}=\frac{1}{n+1}$ when $k\leq n$
and $a_{n,k}=0$ otherwise.
\end{theorem}

We shall write $I_{1}\ll I_{2}$ if there exists a positive constant $K$,
sometimes depended on some parameters, such that $I_{1}\leq KI_{2}$.

\section{Statement of the results}

We start with two propositions.

\begin{proposition}
If $f\in \Omega _{\alpha ,p}$, $n=O\left( r_{n}\right) $ and $q>0,$ then%
\begin{equation*}
\left\Vert \left\{ \frac{1}{r_{n}}\sum_{k=n-r_{n}}^{n-1}\left\vert S_{\frac{%
\alpha k}{2}}f-f\right\vert ^{q}\right\} ^{1/q}\right\Vert _{S^{p}}\ll
\left\Vert f\right\Vert _{S^{p}},
\end{equation*}%
for $n=1,2,...$
\end{proposition}

\begin{proposition}
If $f\in \Omega _{\alpha ,p}$, $n=O\left( r_{n}\right) $ and $q>0,$ then%
\begin{equation*}
\left\Vert \left\{ \frac{1}{r_{n}}\sum_{k=n-r_{n}}^{n-1}\left\vert S_{\frac{%
\alpha k}{2}}f-f\right\vert ^{q}\right\} ^{1/q}\right\Vert _{S^{p}}\ll E_{%
\frac{\alpha \left( n-r_{n}\right) }{2}}(f)_{S^{p}},
\end{equation*}%
for $n=1,2,...$
\end{proposition}

In the special case $p=\infty $ and $f\in C_{2\pi }$\ Proposition 2 reduce
to Lemma 1.

Our main results are following

\begin{theorem}
If $f\in \Omega _{\alpha ,p},$ $p\geq q$, $\left( a_{n,k}\right)
_{k=0}^{\infty }\in GM\left( _{2}\beta \right) $ for all $n$, (\ref{T1}) and 
$\lim_{n\rightarrow \infty }a_{n,0}=0$ hold, then%
\begin{equation}
\left\Vert T_{n,A,\gamma }^{q}f\right\Vert _{S^{p}}\ll \left\{
\sum_{k=0}^{\infty }a_{n,k}E_{\frac{\alpha k}{2^{\left[ c\right] }+1}%
}^{q}(f)_{S^{p}}\right\} ^{1/q},  \label{8}
\end{equation}%
for some $c>1$ and $n=0,1,2,...$, where\ $\gamma =\left( \gamma _{k}\right) $
is a sequence with $\gamma _{k}=\frac{\alpha k}{2}.$
\end{theorem}

\begin{theorem}
If $f\in \Omega _{\alpha ,p},$ $p\geq q$, $\left( a_{n,k}\right)
_{k=0}^{\infty }\in GM\left( _{2}\beta \right) $ for all $n$, (\ref{T1}) and 
$\lim_{n\rightarrow \infty }a_{n,0}=0$ hold, then%
\begin{equation*}
\left\Vert T_{n,A,\gamma }^{q}f\right\Vert _{S^{p}}\ll \left\{
\sum_{k=0}^{\infty }a_{n,k}\omega ^{q}f\left( \frac{\pi }{k+1}\right)
_{S^{p}}\right\} ^{1/q},
\end{equation*}%
for $n=0,1,2,...$, where\ $\gamma =\left( \gamma _{k}\right) $ is a sequence
with $\gamma _{k}=\frac{\alpha k}{2}$.
\end{theorem}

\begin{theorem}
If we additionally suppose that $\left( a_{n,k}\right) _{k=0}^{\infty }\in
MS $ then 
\begin{equation*}
\left\Vert T_{n,A,\gamma }^{q}f\right\Vert _{S^{p}}\ll \left\{
\sum_{k=0}^{\infty }a_{n,k}E_{k}^{p}(f)_{S^{p}}\right\} ^{1/q},
\end{equation*}%
for $n=0,1,2,...$, where\ $\gamma =\left( \gamma _{k}\right) $ is a sequence
with $\gamma _{k}=\frac{\alpha k}{2}$.
\end{theorem}

\begin{remark}
Taking $a_{n,k}=\frac{1}{n+1}$ when $k\leq n$ and $a_{n,k}=0$ otherwise, in
the case $p=\infty $ we obtain the better estimate than this one from \cite[%
Theorem 5]{PT}.
\end{remark}

\section{Proofs of the results}

\subsection{Proof of Proposition 1}

Denote by $S_{k}^{\ast }f$ the sums of the form 
\begin{equation*}
S_{\frac{\alpha k}{2}}f\left( x\right) =\sum_{\left\vert \lambda _{\nu
}\right\vert \leq \frac{\alpha k}{2}}A_{\nu }\left( f\right) e^{i\lambda
_{\nu }x}
\end{equation*}%
\ such that the interval $\left( \frac{\alpha k}{2},\frac{\alpha \left(
k+1\right) }{2}\right) $ does not contain any $\lambda _{\nu }.$ Applying
Lemma 1.10.2 of \cite{BML} we easily verify that%
\begin{equation*}
S_{k}^{\ast }f\left( x\right) -f\left( x\right) =\int_{0}^{\infty }\varphi
_{x}\left( t\right) \Psi _{k}\left( t\right) dt,
\end{equation*}%
where\ $\varphi _{x}\left( t\right) :=f\left( x+t\right) +f\left( x-t\right)
-2f\left( x\right) $ and $\Psi _{k}\left( t\right) =\Psi _{\frac{\alpha k}{2}%
,\frac{\alpha \left( k+1\right) }{2}}\left( t\right) $ i.e.,%
\begin{equation*}
\Psi _{k}\left( t\right) =\frac{4\sin \frac{\alpha t}{4}\sin \frac{\alpha
\left( 2k+1\right) t}{4}}{\pi \alpha t^{2}}
\end{equation*}%
(see also \cite{ASB}, p.41). Evidently, if the interval $\left( \frac{\alpha
k}{2},\frac{\alpha \left( k+1\right) }{2}\right) $ contains a Fourier
exponent $\lambda _{\nu },$ then%
\begin{equation*}
S_{\frac{\alpha k}{2}}f\left( x\right) =S_{k+1}^{\ast }f\left( x\right)
-\left( A_{\nu }\left( f\right) e^{i\lambda _{\nu }x}+A_{-\nu }\left(
f\right) e^{-i\lambda _{\nu }x}\right) .
\end{equation*}%
Since%
\begin{equation*}
\left\{ \sum_{\nu =-\infty }^{\infty }\left\vert A_{\nu }\left( f\right)
\right\vert ^{q}\right\} ^{1/q}\leq \left\Vert f\right\Vert _{B^{p}}\text{ \
for }1<p\leq 2\text{ and }q=\frac{p}{p-1}\text{\ (\cite[p.78]{ABI})}
\end{equation*}%
and%
\begin{equation*}
\left\Vert f\right\Vert _{B^{p}}\leq \left\Vert f\right\Vert _{S^{p}}\text{
\ for }p\geq 1\text{\ (\cite[p. 7]{ADB}),}
\end{equation*}%
where $\left\Vert \cdot \right\Vert _{B^{p}}$\ is the Besicivitch norm,\ we
have%
\begin{equation*}
\left\vert A_{\pm \nu }\left( f\right) \right\vert \leq \left\Vert
f\right\Vert _{S^{p}}\text{ \ \ for }p>1,
\end{equation*}%
whence the deviation 
\begin{equation*}
\frac{1}{r_{n}}\sum_{k=n-r_{n}}^{n-1}\left\vert S_{\frac{\alpha k}{2}%
}f\left( x\right) -f\left( x\right) \right\vert ^{q}
\end{equation*}%
can be estimated from above by%
\begin{equation*}
\frac{1}{r_{n}}\sum_{k=n-r_{n}}^{n-1}\left\vert \int_{0}^{\infty }\varphi
_{x}\left( t\right) \Psi _{k+\kappa }\left( t\right) dt\right\vert ^{q}+%
\frac{1}{r_{n}}\sum_{k=n-r_{n}}^{n-1}\left( \left\Vert f\right\Vert
_{S^{p}}\right) ^{q},
\end{equation*}%
where $\kappa $ equals $0$ or $1.$ Putting $h=2\pi /\left( \alpha n\right) $
we obtain%
\begin{eqnarray*}
\int_{0}^{\infty }\varphi _{x}\left( t\right) \Psi _{k+\kappa }\left(
t\right) dt &=&\left( \int_{0}^{h}+\int_{h}^{nh}+\int_{nh}^{\infty }\right)
\varphi _{x}\left( t\right) \Psi _{k+\kappa }\left( t\right) dt \\
&=&I_{1}\left( k\right) +I_{2}\left( k\right) +I_{3}\left( k\right) .
\end{eqnarray*}%
By elementary calculations we get%
\begin{equation*}
\left\vert I_{1}\left( k\right) \right\vert \leq \frac{\left( 2k+3\right)
\alpha }{4\pi }\int_{0}^{h}\left\vert \varphi _{x}\left( t\right)
\right\vert dt\ll \frac{1}{h}\int_{0}^{h}\left\vert \varphi _{x}\left(
t\right) \right\vert dt
\end{equation*}%
and%
\begin{equation*}
\left\vert I_{3}\left( k\right) \right\vert \leq \int_{nh}^{\infty
}\left\vert \varphi _{x}\left( t\right) \Psi _{k+\kappa }\left( t\right)
\right\vert dt\ll \int_{nh}^{\infty }\frac{\left\vert \varphi _{x}\left(
t\right) \right\vert }{t^{2}}dt.
\end{equation*}%
Therefore%
\begin{equation*}
\frac{1}{r_{n}}\sum_{k=n-r_{n}}^{n-1}\left[ \left\vert I_{1}\left( k\right)
\right\vert +\left\vert I_{3}\left( k\right) \right\vert \right] ^{q}\ll %
\left[ \frac{1}{h}\int_{0}^{h}\left\vert \varphi _{x}\left( t\right)
\right\vert dt+\int_{nh}^{\infty }\frac{\left\vert \varphi _{x}\left(
t\right) \right\vert }{t^{2}}dt\right] ^{q}.
\end{equation*}%
Consequently, we have to estimate the quantity $\frac{1}{r_{n}}%
\sum_{k=n-r_{n}}^{n-1}\left\vert I_{2}\left( k\right) \right\vert ^{q}$. The
inequality of Hausdorf-Young \cite[Chap. XII, Th. 3.3 II]{Z} yields (cf. 
\cite[p. 102]{PT})%
\begin{equation*}
\frac{1}{r_{n}}\sum_{k=n-r_{n}}^{n-1}\left\vert I_{2}\left( k\right)
\right\vert ^{q}\ll \frac{1}{n}\sum_{k=1}^{n}\left\vert I_{2}\left( k\right)
\right\vert ^{q}\ll \frac{1}{n}\left[ \int_{h}^{nh}\frac{\left\vert \varphi
_{x}\left( t\right) \right\vert ^{q^{\prime }}}{t^{q^{\prime }}}dt\right]
^{q/q^{\prime }},
\end{equation*}%
where $q^{\prime }=\frac{q}{q-1}$ and $q\geq 2$.

By monotonicity of $\left\{ \frac{1}{r_{n}}\sum_{k=n-r_{n}}^{n-1}\left\vert
S_{\frac{\alpha k}{2}}f-f\right\vert ^{v}\right\} ^{1/v}$ with respect to $%
v>0,$ 
\begin{eqnarray*}
&&\left\Vert \left\{ \frac{1}{r_{n}}\sum_{k=n-r_{n}}^{n-1}\left\vert S_{%
\frac{\alpha k}{2}}f-f\right\vert ^{v}\right\} ^{1/v}\right\Vert
_{S^{p}}\leq \left\Vert \left\{ \frac{1}{r_{n}}\sum_{k=n-r_{n}}^{n-1}\left%
\vert S_{\frac{\alpha k}{2}}f-f\right\vert ^{q}\right\} ^{1/q}\right\Vert
_{S^{p}} \\
&\ll &\frac{1}{h}\int_{0}^{h}\left\Vert \varphi _{\cdot }\left( t\right)
\right\Vert _{S^{p}}dt+\int_{nh}^{\infty }\frac{\left\Vert \varphi _{\cdot
}\left( t\right) \right\Vert _{S^{p}}}{t^{2}}dt+\left\{ \frac{1}{n}%
\int_{h}^{nh}\frac{\left\Vert \varphi _{\cdot }\left( t\right) \right\Vert
_{S^{p}}^{q^{\prime }}}{t^{q^{\prime }}}dt\right\} ^{1/q^{\prime
}}+\left\Vert f\right\Vert _{S^{p}} \\
&\ll &\frac{1}{h}\int_{0}^{h}\left\Vert f\right\Vert
_{S^{p}}dt+\int_{nh}^{\infty }\frac{\left\Vert f\right\Vert _{S^{p}}}{t^{2}}%
dt+\left\{ \frac{1}{n}\int_{h}^{nh}\frac{\left\Vert f\right\Vert
_{S^{p}}^{q^{\prime }}}{t^{q^{\prime }}}dt\right\} ^{1/q^{\prime
}}+\left\Vert f\right\Vert _{S^{p}} \\
&\ll &\left\Vert f\right\Vert _{S^{p}}\left[ 2+\int_{nh}^{\infty }\frac{1}{%
t^{2}}dt+\left( \frac{1}{n}\int_{h}^{nh}\frac{1}{t^{q^{\prime }}}dt\right)
^{1/q^{\prime }}\right] \ll \left\Vert f\right\Vert _{S^{p}},
\end{eqnarray*}%
for any $v\in \left( 0,q\right] $ such that $q^{\prime }\leq p$.

Thus the desired result follows.$\square $

\subsection{Proof of Proposition 2}

The proof is standard and the estimate follows from that of Proposition 1.
Namely, taking $g_{\sigma }\in B_{\sigma }$ such that $E_{\sigma
}(f)_{S^{p}}=\left\Vert f-g_{\sigma }\right\Vert _{S^{p}}$ we obtain

\begin{eqnarray*}
&&\left\Vert \left\{ \frac{1}{r_{n}}\sum_{k=n-r_{n}}^{n-1}\left\vert S_{%
\frac{\alpha k}{2}}f-f\right\vert ^{q}\right\} ^{1/q}\right\Vert _{S^{p}} \\
&=&\left\Vert \left\{ \frac{1}{r_{n}}\sum_{k=n-r_{n}}^{n-1}\left\vert S_{%
\frac{\alpha k}{2}}f-g_{\sigma }-\left( f-g_{\sigma }\right) \right\vert
^{q}\right\} ^{1/q}\right\Vert _{S^{p}} \\
&=&\left\Vert \left\{ \frac{1}{r_{n}}\sum_{k=n-r_{n}}^{n-1}\left\vert S_{%
\frac{\alpha k}{2}}\left( f-g_{\sigma }\right) -\left( f-g_{\sigma }\right)
\right\vert ^{q}\right\} ^{1/q}\right\Vert _{S^{p}} \\
&\ll &\left\Vert f-g_{\sigma }\right\Vert _{S^{p}},
\end{eqnarray*}%
with $\sigma =\frac{\alpha \left( n-r_{n}\right) }{2},$ and thus our result
follows. $\square $

\subsection{Proof of Theorem 4}

Let%
\begin{equation*}
\left\Vert T_{n,A,\gamma }^{q}f\right\Vert _{S^{p}}=\left\Vert \left\{
\sum_{k=0}^{2^{\left[ c\right] }-1}a_{n,k}\left\vert S_{\frac{\alpha k}{2}%
}f-f\right\vert ^{q}+\sum_{k=2^{[c]}}^{\infty }a_{n,k}\left\vert S_{\frac{%
\alpha k}{2}}f-f\right\vert ^{q}\right\} ^{1/q}\right\Vert _{S^{p}}
\end{equation*}%
\begin{equation*}
\leq \left\Vert \left\{ \sum_{k=0}^{2^{\left[ c\right] }-1}a_{n,k}\left\vert
S_{\frac{\alpha k}{2}}f-f\right\vert ^{q}\right\} ^{1/q}\right\Vert
_{S^{p}}+\left\Vert \left\{ \sum_{m=\left[ c\right] }^{\infty
}\sum_{k=2^{m}}^{2^{m+1}-1}a_{n,k}\left\vert S_{\frac{\alpha k}{2}%
}f-f\right\vert ^{q}\right\} ^{1/q}\right\Vert _{S^{p}}=I_{1}+I_{2}.
\end{equation*}%
for some $c>1$. Using Proposition 2 we obtain, for $p\geq q,$%
\begin{eqnarray*}
I_{1} &\leq &\left\Vert \left\{ \sum_{k=0}^{2^{\left[ c\right] }-1}a_{n,k}%
\frac{k/2+1}{k/2+1}\sum\limits_{l=k/2}^{k}\left\vert S_{\frac{\alpha l}{2}%
}f-f\right\vert ^{q}\right\} ^{1/q}\right\Vert _{S^{p}} \\
&\leq &\left\{ 2^{\left[ c\right] }\sum_{k=0}^{2^{\left[ c\right]
}-1}a_{n,k}\left\Vert \left[ \frac{1}{k/2+1}\sum\limits_{l=k/2}^{k}\left%
\vert S_{\frac{\alpha l}{2}}f-f\right\vert ^{q}\right] ^{1/q}\right\Vert
_{S^{p}}^{q}\right\} ^{1/q} \\
&\ll &\left\{ \sum_{k=0}^{2^{\left[ c\right] }-1}a_{n,k}E_{\frac{\alpha k}{4}%
}^{q}(f)_{S^{p}}\right\} ^{1/q}.
\end{eqnarray*}%
By partial summation, our Proposition 2 gives%
\begin{eqnarray*}
I_{2} &=&\left\Vert \left\{ \sum_{m=[c]}^{\infty }\left[
\sum_{k=2^{m}}^{2^{m+1}-2}\left( a_{n,k}-a_{n,k+1}\right)
\sum_{l=2^{m}}^{k}\left\vert S_{\frac{\alpha l}{2}}f-f\right\vert
^{q}\right. \right. \right. \\
&&\left. \left. \left. +a_{n,2^{m+1}-1}\sum_{l=2^{m}}^{2^{m+1}-1}\left\vert
S_{\frac{\alpha l}{2}}f-f\right\vert ^{q}\right] \right\} ^{1/q}\right\Vert
_{S^{p}}
\end{eqnarray*}%
\begin{eqnarray*}
&\leq &\left\{ \sum_{m=[c]}^{\infty }\left[ \sum_{k=2^{m}}^{2^{m+1}-2}\left%
\vert a_{n,k}-a_{n,k+1}\right\vert \left\Vert \left(
\sum_{l=2^{m}}^{k}\left\vert S_{\frac{\alpha l}{2}}f-f\right\vert
^{q}\right) ^{1/q}\right\Vert _{S^{p}}^{q}\right. \right. \\
&&\left. \left. +a_{n,2^{m+1}-1}\left\Vert \left(
\sum_{l=2^{m}}^{2^{m+1}-1}\left\vert S_{\frac{\alpha l}{2}}f-f\right\vert
^{q}\right) ^{1/q}\right\Vert _{S^{p}}^{q}\right] \right\} ^{1/q}
\end{eqnarray*}%
\begin{eqnarray*}
&\ll &\left\{ \sum_{m=[c]}^{\infty }\left[ 2^{m}\sum_{k=2^{m}}^{2^{m+1}-2}%
\left\vert a_{n,k}-a_{n,k+1}\right\vert E_{\frac{\alpha 2^{m}}{2}%
}^{q}(f)_{S^{p}}\right. \right. \\
&&\left. \left. +2^{m}a_{n,2^{m+1}-1}E_{\frac{\alpha 2^{m}}{2}%
}^{q}(f)_{S^{p}}\right] \right\} ^{1/q} \\
&\ll &\left\{ \sum_{m=[c]}^{\infty }2^{m}E_{\frac{\alpha 2^{m}}{2}%
}^{q}(f)_{S^{p}}\left[ \sum_{k=2^{m}}^{2^{m+1}-2}\left\vert
a_{n,k}-a_{n,k+1}\right\vert +a_{n,2^{m+1}-1}\right] \right\} ^{1/q},
\end{eqnarray*}%
for $p\geq q$.

Since (\ref{5}) holds, we have%
\begin{eqnarray*}
&&a_{n,s+1}-a_{n,r} \\
&\leq &\left\vert a_{n,r}-a_{n,s+1}\right\vert \leq \sum_{k=r}^{s}\left\vert
a_{n,k}-a_{n,k+1}\right\vert \\
&\leq &\sum_{k=2^{m}}^{2^{m+1}-2}\left\vert a_{n,k}-a_{n,k+1}\right\vert \ll
\sum\limits_{k=[2^{m}/c]}^{[c2^{m}]}\frac{a_{n,k}}{k}\text{ \ \ }\left(
2\leq 2^{m}\leq r\leq s\leq 2^{m+1}-2\right) ,
\end{eqnarray*}%
whence%
\begin{equation*}
a_{n,s+1}\ll a_{n,r}+\sum\limits_{k=[2^{m}/c]}^{[c2^{m}]}\frac{a_{n,k}}{k}%
\text{ \ }\left( 2\leq 2^{m}\leq r\leq s\leq 2^{m+1}-2\right) .
\end{equation*}%
Consequently,%
\begin{eqnarray*}
2^{m}a_{n,2^{m+1}-1} &=&\frac{2^{m}}{2^{m}-1}%
\sum_{r=2^{m}}^{2^{m+1}-2}a_{n,2^{m+1}-1} \\
&\ll &\sum_{r=2^{m}}^{2^{m+1}-2}\left(
a_{n,r}+\sum\limits_{k=[2^{m}/c]}^{[c2^{m}]}\frac{a_{n,k}}{k}\right) \\
&\ll
&\sum_{r=2^{m}}^{2^{m+1}-1}a_{n,r}+2^{m}\sum\limits_{k=[2^{m}/c]}^{[c2^{m}]}%
\frac{a_{n,k}}{k}
\end{eqnarray*}%
and therefore%
\begin{equation*}
I_{2}\ll \left\{ \sum_{m=[c]}^{\infty }\left[ 2^{m}E_{\frac{\alpha 2^{m}}{2}%
}^{q}(f)_{S^{p}}\sum\limits_{k=[2^{m}/c]}^{[c2^{m}]}\frac{a_{n,k}}{k}+E_{%
\frac{\alpha 2^{m}}{2}}^{q}(f)_{S^{p}}\sum_{k=2^{m}}^{2^{m+1}-1}a_{n,k}%
\right] \right\} ^{1/q}.
\end{equation*}

Finally, by elementary calculations we get%
\begin{eqnarray*}
I_{2} &\ll &\left\{ \sum_{m=[c]}^{\infty }\left[ 2^{m}E_{\frac{\alpha 2^{m}}{%
2}}^{q}(f)_{S^{p}}\sum\limits_{k=2^{m-\left[ c\right] }}^{2^{m+\left[ c%
\right] }}\frac{a_{n,k}}{k}+E_{\frac{\alpha 2^{m}}{2}}^{q}(f)_{S^{p}}%
\sum_{k=2^{m}}^{2^{m+1}}a_{n,k}\right] \right\} ^{1/q} \\
&\ll &\left\{ \sum_{m=[c]}^{\infty }E_{\frac{\alpha 2^{m}}{2}%
}^{q}(f)_{S^{p}}\sum\limits_{k=2^{m-\left[ c\right] }}^{2^{m+\left[ c\right]
}}a_{n,k}\right\} ^{1/q} \\
&=&\left\{ \sum_{m=[c]}^{\infty }E_{\frac{\alpha 2^{m}}{2}%
}^{q}(f)_{S^{p}}\sum\limits_{k=2^{m-\left[ c\right] }}^{2^{m}-1}a_{n,k}+%
\sum_{m=[c]}^{\infty }E_{\frac{\alpha 2^{m}}{2}}^{q}(f)_{S^{p}}\sum%
\limits_{k=2^{m}}^{2^{m+\left[ c\right] }}a_{n,k}\right\} ^{1/q}
\end{eqnarray*}%
\begin{equation*}
\ll \left\{ \sum_{m=[c]}^{\infty }\sum\limits_{k=2^{m-\left[ c\right]
}}^{2^{m}-1}a_{n,k}E_{\frac{\alpha k}{2}}^{q}(f)_{S^{p}}+\sum_{m=[c]}^{%
\infty }\sum\limits_{k=2^{m}}^{2^{m+\left[ c\right] }}a_{n,k}E_{\frac{%
\alpha k}{2^{\left[ c\right] +1}}}^{q}(f)_{S^{p}}\right\} ^{1/q}
\end{equation*}%
\begin{equation*}
\ll \left\{ \sum_{m=[c]}^{\infty }\sum\limits_{k=2^{m-\left[ c\right]
}}^{2^{m}-1}a_{n,k}E_{\frac{\alpha k}{2}}^{q}(f)_{S^{p}}+\sum_{m=[c]}^{%
\infty }\sum\limits_{k=2^{m}}^{2^{m+\left[ c\right] }-1}a_{n,k}E_{\frac{%
\alpha k}{2^{\left[ c\right] +1}}}^{q}(f)_{S^{p}}+\sum_{m=[c]}^{\infty }E_{%
\frac{\alpha 2^{m}}{2}}^{q}(f)_{S^{p}}a_{n,2^{m+\left[ c\right] }}\right\}
^{1/q}
\end{equation*}%
\begin{eqnarray*}
&=&\left\{ \sum_{m=[c]}^{\infty }\sum\limits_{r=1}^{\left[ c\right]
}\sum\limits_{k=2^{m-r}}^{2^{m-r+1}-1}a_{n,k}E_{\frac{\alpha k}{2}%
}^{q}(f)_{S^{p}}+\sum_{m=[c]}^{\infty }\sum\limits_{r=0}^{\left[ c\right]
-1}\sum\limits_{k=2^{m+r}}^{2^{m+r+1}-1}a_{n,k}E_{\frac{\alpha k}{2^{\left[
c\right] +1}}}^{q}(f)_{S^{p}}\right. \\
&&+\left. \sum_{m=[c]}^{\infty }E_{\frac{\alpha 2^{m}}{2}%
}^{q}(f)_{S^{p}}a_{n,2^{m+\left[ c\right] }}\right\} ^{1/q}
\end{eqnarray*}%
\begin{eqnarray*}
&\leq &\left\{ \sum\limits_{r=1}^{\left[ c\right] }\sum\limits_{k=2^{\left[
c\right] -r}}^{\infty }a_{n,k}E_{\frac{\alpha k}{2}}^{q}(f)_{S^{p}}+\sum%
\limits_{r=0}^{\left[ c\right] -1}\sum\limits_{k=2^{\left[ c\right]
+r}}^{\infty }a_{n,k}E_{\frac{\alpha k}{2^{\left[ c\right] +1}}%
}^{q}(f)_{S^{p}}+\sum\limits_{k=2^{2\left[ c\right] }}^{\infty }a_{n,k}E_{%
\frac{\alpha 2^{m}}{2}}^{q}(f)_{S^{p}}\right\} ^{1/q} \\
&\ll &\left\{ \sum\limits_{k=0}^{\infty }a_{n,k}E_{\frac{\alpha k}{2^{\left[
c\right] +1}}}^{q}(f)_{S^{p}}\right\} ^{1/q}.
\end{eqnarray*}

Thus we obtain the desired result. $\square $

\subsection{Proof of Theorem 5}

The proof follows by the Jackson type theorem%
\begin{equation*}
E_{\sigma }(f)_{S^{p}}\ll \omega f\left( \frac{1}{\sigma }\right) _{S^{p}}
\end{equation*}%
and basic properties of the modulus of continuity $\omega f\left( \cdot
\right) _{S^{p}}$ .$\square $

\subsection{Proof of Theorem 6}

If $\left( a_{n,k}\right) _{k=0}^{\infty }\in MS$ then $\left(
a_{n,k}\right) _{k=0}^{\infty }\in GM\left( _{2}\beta \right) $ and using
Theorem 4 we obtain%
\begin{eqnarray*}
\left\Vert T_{n,A,\gamma }^{q}f\right\Vert &\ll &\left\{
\sum\limits_{k=0}^{\infty }a_{n,k}E_{\frac{\alpha k}{2^{\left[ c\right] +1}}%
}^{q}(f)_{S^{p}}\right\} ^{1/q}=\left\{ \sum_{k=0}^{\infty
}\sum\limits_{m=k2^{\left[ c\right] }}^{\left( k+1\right) 2^{\left[ c\right]
}-1}a_{n,m}E_{\frac{\alpha k}{2^{\left[ c\right] +1}}}^{q}(f)_{S^{p}}\right%
\} ^{1/q} \\
&\leq &\left\{ \sum_{k=0}^{\infty }E_{\frac{\alpha k}{2}}^{q}(f)_{S^{p}}%
\sum\limits_{m=k2^{\left[ c\right] }}^{\left( k+1\right) 2^{\left[ c\right]
}-1}a_{n,m}\right\} ^{1/q}\leq \left\{ \sum_{k=0}^{\infty }2^{\left[ c\right]
}E_{\frac{\alpha k}{2}}^{q}(f)_{S^{p}}a_{n,k2^{\left[ c\right] }}\right\}
^{1/q} \\
&\leq &\left\{ 2^{\left[ c\right] }\sum_{k=0}^{\infty }E_{\frac{\alpha k}{2}%
}^{q}(f)_{S^{p}}a_{n,k}\right\} ^{1/q}\ll \left\{ \sum_{k=0}^{\infty }E_{%
\frac{\alpha k}{2}}^{q}(f)_{S^{p}}a_{n,k}\right\} ^{1/q}.
\end{eqnarray*}%
This ends of our proof. $\square $

\end{document}